\documentclass[11pt]{article}
\usepackage{amssymb,amsmath,bbm,latexsym,amscd}
\newtheorem{theorem}[equation]{Theorem}

\newtheorem{lemma}[equation]{Lemma}
\newtheorem{proposition}[equation]{Proposition}

\newtheorem{remark}[equation]{Remark}
\numberwithin{equation}{section}

\newcommand{\ot}{\otimes}
\newcommand{\ol}{\overline}

\newcommand{\Z}{\mathbb{Z}}

\newcommand{\fg}{{\mathfrak g}}

\newcommand{\Rep}{\hbox{Rep}\,}
\newcommand{\Max}{\hbox{Max}}

\newcommand{\gs}{\sigma}

\newcommand{\proof}{{\bf Proof\ \ }}
\newcommand{\qed}{\hfill $\Box$}

\newcommand{\End}{\hbox{End}}

\newcommand{\cL}{{\mathcal L}}

\newcommand{\Aut}{{\rm Aut}}

\newcommand\Spec{\text{\rm Spec}\,}

\newcommand{\cI}{\mathcal{I}}
\newcommand{\cM}{\mathcal{M}}
\newcommand{\ev}{\mathrm{ev}}
\newcommand{\supp}{\hbox{supp}\,}

\title{Representations of Twisted Current Algebras}

\date{}

\author{Michael Lau\thanks{Funding from the Natural Sciences and
Engineering Research Council of Canada is gratefully acknowledged.} \vspace{0.3cm}
\\{\small Universit\'e Laval}\\
{\small D\'epartement de math\'ematiques et de statistique}\\ {\small Qu\'ebec, QC, Canada G1V 0A6}\\ {\small Email:\ Michael.Lau@mat.ulaval.ca}\vspace{0.1cm}}

\begin{document}
\maketitle

\begin{small}
\noindent {\bf Abstract.} We use evaluation representations to give a complete classification of the finite-dimensional simple modules of twisted current algebras.  This generalizes and unifies recent work on multiloop algebras, current algebras, equivariant map algebras, and twisted forms.

\bigskip




\end{small}
\maketitle

\section{Introduction}  Evaluation representations were introduced in the late 1980s to classify representations of loop algebras \cite{chari86, cp88,rao93}.  This powerful technique was later adapted to study representations of multiloop algebras \cite{rao01,batra,multiloop} and various other generalizations and related Lie algebras.  (See \cite{feigenloktev04,charifourierkhandai10,hartwig07,dateroan00}, for instance.)  Two very recent papers use evaluation representations to classify the finite-dimensional simple modules of {\em equivariant map algebras} and {\em twisted forms} \cite{NSS,repforms}.  While these two families contains all the previously mentioned examples in the literature, twisted forms also include {\em Margaux algebras} (which are not equivariant map algebras) and equivariant map algebras include {\em Onsager algebras} (which are not twisted forms).  Equivariant map algebras and twisted forms are both examples of {\em twisted current algebras}, fixed point subalgebras of tensor products of Lie algebras and associative algebras under a finite group action.  The goal of this short paper is to classify the finite-dimensional simple modules of twisted current algebras up to isomorphism.  This gives a unified perspective on evaluation representations, recovering the main results of \cite{NSS} and \cite{repforms} as special cases.

The work in \cite{repforms} suggests that the common theme in the evaluation module approach is not the equivariant maps studied in \cite{multiloop,NSS}, but rather the {\em invariants} of current algebras under a finite group action.  Let $\fg$ be a finite-dimensional simple Lie algebra over an algebraically closed field $k$ of characteristic zero, and let $X=\Spec S$ be a reduced affine $k
$-scheme of finite type with coordinate algebra $S$.  The $\fg$-valued regular functions on $X$ form a {\em current Lie algebra} $\fg\ot_k S$ over $k$, with pointwise Lie bracket $[x\ot r,y\ot s]=[x,y]\ot rs$, for $x,y\in\fg$ and $r,s\in S$.  When a finite group $\Gamma$ acts by $k$-Lie algebra automorphisms on $\fg\ot_k S$, the {\em twisted current algebra} is simply the fixed point subalgebra $\cL=(\fg\ot_k S)^\Gamma$ under the action of $\Gamma$.

Special cases include both the equivariant map algebras considered in \cite{NSS} and the twisted forms of algebras studied in \cite{repforms}.  Indeed, we obtain the equivariant map algebras when $\Gamma$ acts by automorphisms on both the Lie algebra $\fg$ and the $k$-algebra $S$, such that $${}^\gamma(x\ot s)={}^\gamma x\ot {}^{\gamma}s\ (x\in\fg,\ s\in S,\ \gamma\in \Gamma).$$
The twisted forms of \cite{repforms} appear when $\Gamma$ acts on $\fg\ot S$ via a group action on $S$ and a $1$-cocycle $u:\ \Gamma\rightarrow \Aut_{S-Lie}(\fg\ot_k S)$, such that 
$${}^\gamma(x\ot s)=u_\gamma(x\ot{}^\gamma s),$$
and the ring extension of $S$ over $S^\Gamma=\{s\in S\ :\ {}^\gamma s=s\mbox{\ for all\ }\gamma\in\Gamma\}$ is Galois.  

The classification of finite-dimensional simple modules of a twisted current algebra $\cL=(\fg\ot_k S)^\Gamma$ begins with the observation that any such representation $\phi:\ \cL\rightarrow \End\,V$ descends to a representation of $\cL/\ker\phi$, a finite-dimensional reductive Lie algebra.  Up to twisting by a character, the quotient map $\cL\rightarrow\cL/\ker\phi\subseteq \mbox{End}\,V$ factors through the image of an evaluation homomorphism
$$\ev_{\underline{M}}:\ \cL\hookrightarrow\fg\ot_k S\rightarrow\left(\fg\ot_k S/M_1
\right)\oplus\cdots\oplus\left(\fg\ot_k S/M_r\right)\cong\fg^{\oplus r},$$
for some family of maximal ideals $\underline{M}=\{M_1,\ldots,M_r\}\subseteq\Max\, S$.  The isomorphism class of the representation $(V,\phi)$ is then determined (up to twisting by a character) by the choice of $\underline{M}$ and by the isomorphism class of $V$, viewed as a finite-dimensional simple $\ev_{\underline{M}}(\cL)$-module.

In the case of equivariant map algebras, the image of the evaluation homomorphism is a direct sum of isotropy subalgebras $\fg^{M_i}$ with respect to an action of subgroups $\Gamma^{M_i}\subseteq\Gamma$ fixing maximal ideals $M_i$.  A priori, the definition of a twisted current algebra $\cL=(\fg\ot_k S)^\Gamma$ contains neither an action of $\Gamma$ on $\fg$ nor on $S$.  However, the action of $\Gamma$ on $\fg\ot_k S$ induces an action on the maximal ideals of $S$, and thus on $S$ itself.  This lets us define a $1$-cocycle $u:\ \Gamma\rightarrow\Aut_{S-Lie}(\fg\ot_k S)$, so that the action of $\Gamma$ on $\fg\ot_k S$ can be described as in the case of twisted forms: ${}^\gamma(x\ot s)=u_\gamma(x\ot{}^\gamma s).$  While the general ideas of the evaluation module approach remain the same, the fact that the extension $S/S^\Gamma$ need not be Galois (and thus $\cL\ot_{S^\Gamma} S$ need not be isomorphic to $\fg\ot S$) invalidates some of the most important descent-based arguments in this context.

Instead, we adopt the approach presented in \cite{NSS}.  The major obstacle here is the reliance on an action of $\Gamma$ on the finite-dimensional Lie algebra $\fg$ to define equivariant map algebras.  Such an action no longer exists in the context of twisted current algebras.  We finesse this difficulty by using {\em local twisted actions}, an idea motivated by ideas from affine group schemes.  If $\Gamma\subseteq \Aut_{S-Lie}(\fg\ot_k S)$, then $\Gamma$ is an affine group scheme over $S$, with $\Gamma(A)\subseteq\Aut_{A-Lie}(\fg(A))$ given by the $A$-linear extension of the $\Gamma$-action on $\fg\ot_k S$ to $\fg(A)=(\fg\ot_k S)\ot_S A$.  In particular, if $A=S/M$ for some maximal ideal $M$ of $S$, then $\fg(A)\cong\fg$, and we have a group homomorphism $\ev_M$ associated to each $M\in\Max\, S$:
\begin{eqnarray*}
\Gamma&\stackrel{\ev_M}{\longrightarrow}&\Gamma(S/M)\subseteq\Aut_{k-Lie}(\fg)\\
\alpha&\longmapsto&\alpha(M).
\end{eqnarray*}

When working with twisted current algebras, the action of $\Gamma$ is only $k$-linear, in general, and each evaluation map $\ev_M$ gives a twisted action on $\fg$:
$$(\alpha\beta)(M)=\alpha(M) {}^\alpha\beta(M)\mbox{\quad\quad\quad\quad}(\alpha,\beta\in\Aut_{k-Lie}(\fg\ot_k S),\ M\in\Max S).$$
See Section 3 for details.  This local twisted action turns out to be the crucial ingredient to adapt the proofs of \cite{NSS} to the framework of twisted forms, and more generally, to twisted current algebras.

The methods in this paper can also be used to study the abelian tensor category $\mathcal{C}$ of finite-dimensional modules of twisted current algebras, a topic we plan to explore in a future paper \cite{fdreps}.  In particular, this will classify the blocks of $\mathcal{C}$, and consequently the blocks of twisted forms of algebras.








\section{Twisted current algebras}\label{ketto}
Throughout this paper, $k$ will denote an algebraically closed field of characteristic zero, $\fg$ will be a finite-dimensional simple Lie algebra over $k$, $X$ will be a reduced affine $k$-scheme of finite type, and the {\em (untwisted) current algebra} $\mathfrak{M}=\mathfrak{M}(X,\fg)$ will be the space of $\fg$-valued regular functions on $X$.  This is a Lie algebra under the pointwise Lie bracket $[f,g](M)=[f(M),g(M)],$ for all $f,g\in\mathfrak{M}$ and $k$-rational points $M\in X$. Let $\Gamma$ be a finite group acting on $\mathfrak{M}$ by $k$-Lie algebra automorphisms.  The subalgebra of $\Gamma$-invariants,
$$\mathfrak{M}(X,\fg)^\Gamma=\{f\in\mathfrak{M}\ :\ {}^\gamma f=f\hbox{\ for all\ } \gamma\in\Gamma\},$$
is called the {\em twisted current algebra} associated with $\mathfrak{M}$ and $\Gamma$.  

If $S$ is the coordinate algebra of $X$, we denote by $s(M)\in k$ the residue class of $s\in S$ modulo any maximal ideal ($k$-rational point) $M\in\Max\,S$:
$$s(M)+M=s+M\in S/M\cong k.$$
The twisted current algebra $\mathfrak{M}(X,\fg)^\Gamma$ is isomorphic to, and will henceforth be identified with, the Lie algebra of fixed points
$$\cL=\left(\fg\ot_k S\right)^\Gamma=\{z\in\fg\ot_k S\ :\ {}^\gamma z=z\},$$
where the action of $\Gamma$ on $\fg\ot_k S$ is obtained by identifying $\fg\ot_k S$ with $\mathfrak{M}(X,\fg)$ via
\begin{equation}\label{evaluation}
\left(\sum x_i\ot_k s_i\right)(M):=\sum s_i(M)x_i,
\end{equation}
for all $x_i\in\fg$, $s_i\in S$, and $k$-rational points $M$ in $X$.  Twisted current algebras are then precisely the $\Gamma$-fixed subalgebras of $k$-Lie algebras $\fg\ot_k S$, where $S$ is a unital, commutative, associative, and reduced $k$-algebra of finite type.  However, note that while the untwisted current algebra $\fg \ot_k S$ is naturally a Lie algebra over the base ring $S$, the action of $\Gamma$, being only $k$-linear, defines only a $k$-Lie algebra structure on the twisted current algebra $\cL=\left(\fg\ot_k S\right)^\Gamma$.

\begin{proposition}\label{1-cocycle}  Let $\cL=\left(\fg\ot_k S\right)^\Gamma$ be the twisted current algebra associated to the action of $\Gamma$ on $\fg\ot_k S$.  Then there is a group action of $\Gamma$ on $S$ and a crossed homomorphism $u:\ \Gamma\rightarrow\Aut_{S-Lie}(\fg\ot_k S)$ such that
$${}^\gamma(x\ot_k s)=u_\gamma(x\ot_k{}^\gamma s),$$
for all $\gamma\in\Gamma,$ $x\in\fg$, and $s\in S$.
\end{proposition}
\proof Let $\gamma\in\Gamma$, and let $P$ be an ideal of $S$.  Since $\Gamma$ acts by automorphisms on $\fg\ot_k S$, the image ${}^\gamma(\fg\ot_k P)$ of $(\fg\ot_k P)$ under the action of $\gamma$ is an ideal of $\fg\ot_k S$.  Every such ideal is of the form $\fg\ot_k I$ for some ideal $I$ of $S$, so 
$${}^\gamma(\fg\ot_k P)=\fg\ot_k{}^\gamma P,$$
for some ideal ${}^\gamma P\subseteq S$.  Restricted to the set $\Max\,S$ of maximal ideals of $S$, this gives a group action of $\Gamma$ on $\Max\,S$.

Since $S$ is reduced and of finite type over an algebraically closed field, the Jacobson radical of $S$ is trivial.  The action of $\Gamma$ on $\Max\,S$ thus determines an action of $\Gamma$ on $S$ by $k$-algebra automorphisms.  Explicitly, the image ${}^\gamma s$ of each $s\in S$ under $\gamma\in\Gamma$ is determined by its reduction modulo the maximal ideals $M\in\Max\,S$:
\begin{equation}\label{1star}
{}^\gamma s(M)=s({}^{\gamma^{-1}}M).
\end{equation} 

For each $\gamma\in\Gamma$, let $u_\gamma:\ \fg\ot_k S\rightarrow \fg\ot_k S$ be the $k$-linear map defined by
$$u_\gamma(x\ot s)={}^\gamma(x\ot{}^{\gamma^{-1}}s),$$
for all $x\ot s\in\fg\ot_k S$.  It is straightforward to verify that $u_\gamma$ is a $k$-Lie algebra automorphism of $\fg\ot_k S$.  

We now show that the automorphism $u_\gamma$ is actually $S$-linear.  Consider $r,s\in S$, $x\in\fg$, and $\gamma\in\Gamma$.  For each $P\in\Max\,S$, there exists $t\in{}^{\gamma^{-1}}P$ such that $s=s({}^{\gamma^{-1}}P)+t$.  Therefore,
\begin{align*}
{}^\gamma(x\ot rs)(P)&={}^\gamma\left(x\ot (r)(s({}^{\gamma^{-1}}P)\right)(P)+{}^\gamma(x\ot rt)(P) \\
&=s({}^{\gamma^{-1}}P){}^\gamma(x\ot r)(P)+{}^\gamma(x\ot rt)(P).
\end{align*}
But $rt\in{}^{\gamma^{-1}}P$, so ${}^\gamma(x\ot rt)\in{}^\gamma(\fg\ot{}^{\gamma^{-1}}P)=\fg\ot P$.  Thus $ {}^\gamma(x\ot rt)(P)=0$ and 
$${}^\gamma(x\ot rs)(P)={}^\gamma s(P){}^\gamma(x\ot r)(P).$$
Expanding with respect to a $k$-basis $\{x_i\}$ of $\fg$, we write 
$${}^\gamma(x\ot rs)=\sum_i x_i\ot m_i\mbox{\quad and\quad}{}^\gamma(x\ot r)=\sum_i x_i\ot r_i,$$
for some $m_i,r_i\in S$.  Thus ${}^\gamma(x\ot rs)(P)=\sum_i m_i(P)x_i$ and ${}^\gamma s(P){}^\gamma(x\ot r)(P)={}^\gamma s(P)\sum r_i(P)x_i$, so $m_i(P)={}^\gamma s(P)r_i(P)$ for all $P\in\Max\,S$ and for all $i$.  Since the Jacobson radical of $S$ is zero, $m_i={}^\gamma sr_i$ for all $i$, and
\begin{equation}\label{S-semilinear}
{}^\gamma(x\ot rs)={}^\gamma s{}^\gamma(x\ot r).
\end{equation}
In particular, 
\begin{align}
u_\gamma(x\ot rs)&={}^\gamma\left(x\ot{}^{\gamma^{-1}}(rs)\right)\nonumber\\
&={}^\gamma\left(x\ot{}^{\gamma^{-1}}r{}^{\gamma^{-1}}s\right)\nonumber\\
&=s{}^\gamma(x\ot{}^{\gamma^{-1}}r)\nonumber\\
&=su_\gamma(x\ot r),\label{2stars}
\end{align}
so $u_\gamma$ is $S$-linear.

The finite group $\Gamma$ acts on the group $\Aut_{S-Lie}(\fg\ot_k S)$ of $S$-Lie algebra automorphisms in the usual way:
\begin{equation}\label{action-on-autos}
{}^\gamma\phi=(1\ot\gamma)\circ\phi\circ(1\ot\gamma^{-1}),
\end{equation} 
where $(1\ot\gamma)(x\ot s)=x\ot{}^\gamma s$ for all $\gamma\in\Gamma$, $\phi\in \Aut_{S-Lie}(\fg\ot_k S)$, and $x\ot s\in\fg\ot_k S$.  Therefore,
\begin{align*}
u_\gamma{}^\gamma u_\eta(x\ot s)&=u_\gamma\left((1\ot \gamma){}^\eta(x\ot{}^{\eta^{-1}\gamma^{-1}}s)\right)\\
&={}^\gamma\left((1\ot\gamma^{-1})\circ(1\ot\gamma){}^\eta(x\ot{}^{\eta^{-1}\gamma^{-1}}s)\right)\\
&={}^{\gamma\eta}\left(x\ot{}^{\eta^{-1}\gamma^{-1}}s)\right)\\
&=u_{\gamma\eta}(x\ot s),
\end{align*}
for all $x\ot s\in\fg\ot S$ and $\gamma,\eta\in\Gamma$.  Thus $u_{\gamma\eta}=u_\gamma{}^\gamma u_\eta,$
and $u:\ \Gamma\rightarrow\Aut_{S-Lie}(\fg\ot_k S)$ is a crossed homomorphism.\qed

\bigskip

Let $R=S^\Gamma=\{s\in S\ :\ {}^\gamma s=s\mbox{\ for all\ }\gamma\in\Gamma\}$ be the subalgebra of $\Gamma$-invariants in $S$.  Then $r\cL\subseteq\cL$ for all $r\in R$, so $\cL$ has the structure of a Lie algebra over $R$, and it follows that all maximal $k$-Lie algebra ideals $\mathcal{I}$ of $\cL$ which do not contain the derived subalgebra $[\cL,\cL]$ are actually $R$-Lie algebra ideals of $\cL$.  See \cite[Lemma 4.16]{NSS} or \cite[Lemma 2.7]{repforms} for details.  Note that these are precisely the ideals $\mathcal{I}$ of $\cL$ for which the quotient $\cL/\mathcal{I}$ is a simple Lie algebra.  We denote this set of $R$-ideals by $\Max\,\cL$.  For each $M\in\Max\,S$, we write $\Gamma^M=\{\gamma\in\Gamma\ :\ {}^\gamma M=M\}$ for the isotropy subgroup of $M$, and we denote by $\psi(M):\ \fg\rightarrow\fg$ the map $\psi(M)x=(\psi(x\ot 1))(M)$ for each $\psi\in\Aut_{k-Lie}(\fg\ot S)$ and $x\in\fg$.  The finite-dimensional Lie algebra $\fg^M=\{x\in\fg\ :\ \gamma(M)x=x\mbox{\ for all\ }\gamma\in\Gamma^M\}$ will reappear later in this section as the image of the evaluation map $\ev_M:\ \cL\rightarrow\fg$ sending an element $\alpha\in\cL$ to $\alpha(M)$, as defined in (\ref{evaluation}).  

\begin{lemma}
For each $R$-ideal $\cM$ of $\cL$, let $I(\cM)=\{ r\in S\ :\ r\cL\subseteq \cM\}$.  Then $I(\cM)\in\Max\,R$ for all $\cM\in\Max\,\cL$.
\end{lemma}
\proof Suppose $\cM\in\Max\,\cL$.  Then $\cL/\cM$ is a simple module for the (countable dimensional) Lie algebra $\cL$, and $I(\cM)$ is the kernel of the map $m:\ R\rightarrow \End_{\cL-mod}(\cL/\cM)$, $m:\ r\mapsto m_r$, where $m_r:\ \cL/\cM\rightarrow \cL/\cM$,\   $\overline{x}\mapsto\overline{rx}.$  By Schur's Lemma, we see that each $m_r$ is a scalar multiplication, and $I(\cM)=\ker\,m$ is a maximal ideal of $R$.\qed

\bigskip

\begin{remark} {\em
If the extension $S/R$ is Galois, then the map $I:\ \Max\,\cL\rightarrow\Max\,R$ is a bijection, as follows from \cite[Theorem 2.9]{repforms}.  In the context of twisted current algebras, this is not true in general.  For example, let $\fg=\mathfrak{sl}_2(k)$, $S=k\oplus k\oplus k$, and $\Gamma$ the Klein-IV group $\Z_2\oplus\Z_2$.  Choose Chevalley generators $e,f,h$ for $\mathfrak{sl}_2(k)$.  Let $(a,b)\in\Gamma$ act on $x\ot s$ for $x\in\fg$ and $s=(s_1,s_2,s_3)\in S$ by 
$$(a,b)\cdot(x\ot s)=\gs_1^a\circ\gs_2^b(x)\ot \gamma^{a+b}(s),$$ 
where $\gs_\ell(e)=(-1)^\ell f$, $\gs_\ell(f)=(-1)^\ell e$, and $\gs_\ell(h)=-h$ for $\ell=1,2$, and $\gamma(s_1,s_2,s_3)=(s_2,s_1,s_3)$.  The maximal ideal $J=\{(a,a,0)\ :\ a\in k\}$ of $R=\{(a,a,b)\ :\ a,b\in k\}$ has the property that $J\cL=\cL$, where $\cL=(\fg\ot_k S)^\Gamma=\hbox{Span}_k\,\{ h\ot_k(1,-1,0)\}$, so there is no maximal ideal $\mathcal{M}\subset \cL$ for which $I(\mathcal{M})=J$.
}\end{remark}

We now prove the following proposition using a modification of the argument in \cite[Proposition 5.2]{NSS}.  Because the proposition is nontrivial and plays an important role in what follows, we give a detailed proof in the context of twisted current algebras for completeness.

\begin{proposition}\label{fgM-epi}
Let $\cM_1,\ldots,\cM_s\in\Max\,\cL$ be distinct ideals for which $I(\cM_1)=\cdots=I(\cM_s)$.  Then there is a natural Lie algebra epimorphism $\fg^M\rightarrow\cL/\cM_1\times\cdots\times\cL/\cM_s$ for any $M\in\Max\,S$ lying over $I=I(\cM_1)$.
\end{proposition}
\proof Since $\Gamma$ fixes the elements of $R$, it leaves $I$ and $IS$ (setwise) invariant.  The exact sequence $0\rightarrow IS\rightarrow S\rightarrow S/IS\rightarrow 0$ induces an exact sequence of $\Gamma$-modules:
$$0\rightarrow \fg\ot_k IS\rightarrow \fg\ot_k S\rightarrow \fg\ot_k S/IS\rightarrow 0.$$
Taking $\Gamma$-invariants is an exact functor, so 
\begin{equation}\label{reduction-iso}
\cL/I\cL=(\fg\ot_k S)^\Gamma/(\fg\ot_k IS)^\Gamma\cong\left(\fg\ot_k S/IS\right)^\Gamma.
\end{equation}
By definition, $I\cL\subseteq \cM_1\cap\cdots\cap\cM_s$, so we have a surjective homomorphism
$$\psi:\ (\fg\ot_k S/IS)^\Gamma\rightarrow\cL/(\cM_1\cap\cdots\cap\cM_s).$$

The extension $S/R$ is integral \cite[\S 1.9 Proposition 22]{bourbaki-commalg5}, so $\Gamma$ acts transitively on the set of (finitely many) maximal ideals $M_1,\ldots,M_r\in\Max\,S$ lying over $I\in\Max\,R$, and $\sqrt{IS}=\bigcap_{i=1}^r M_i$ by \cite[\S 2.1 Proposition 1 and \S 2.2 Th\'eor\`eme 2]{bourbaki-commalg5}.  The action of $\Gamma$ on $\fg\ot S/\sqrt{IS}$ thus transitively permutes the summands $\fg\ot_k S/M_i$ in the direct sum
$$\bigoplus_i(\fg\ot_k S/M_i)\cong\fg\ot_k(S/\cap M_i)\cong\fg\ot_k(S/\sqrt{IS}).$$
For any $M=M_i$ lying over $I$, we see that projection onto the $i$th component $\fg\ot_k S/M_i$ induces an isomorphism of $\Gamma$-invariants: $(\fg\ot S/\sqrt{IS})^\Gamma\stackrel{\pi}{\rightarrow}(\fg\ot_k S/M)^{\Gamma^M}$, since taking $\Gamma$-invariants is an exact functor.  

The evaluation map
\begin{align*}
\ev_M:\ (\fg\ot_k S/M)^{\Gamma^M}&\rightarrow\fg\\
\sum x_i\ot (s_i+M)&\mapsto\sum s_i(M)x_i
\end{align*}
has image contained in $\fg^M$.  Indeed, for any $\sum x_i\ot (s_i+M)\in(\fg\ot_k S/M)^\Gamma$ and $\gamma\in\Gamma^M$, we see that $s_i(M)=s_i+m_i$ for some $m_i\in M$ and
\begin{align*}
\gamma(M)\left(\ev_M\left(\sum x_i\ot (s_i+M)\right)\right)&=\gamma(M)\sum s_i(M)x_i\\
&=\sum {}^\gamma(x_i\ot s_i(M))(M)\\
&={}^\gamma\left(\sum x_i\ot(s_i+m_i)\right)(M)\\
&=\ev_M\left(\sum x_i\ot s_i\right)+ {}^\gamma\left(\sum x_i\ot m_i\right)(M).
\end{align*}
Since $\gamma\in\Gamma^M$, we see that ${}^\gamma\left(\sum x_i\ot m_i\right)\in{}^\gamma(\fg\ot_k M)=\fg\ot_k M$, so ${}^\gamma\left(\sum x_i\ot m_i\right)(M)=0$, $\ev_M\left(\sum x_i\ot (s_i+M)\right)$ is fixed under $\gamma(M)$, and the image of $\ev_M$ is contained in $\fg^M$.  It is then straightforward to verify that $\ev_M:\ (\fg\ot_k S/M)^{\Gamma^M}\rightarrow\fg^M$ is a Lie algebra isomorphism.

Let $\alpha=\sum x_i\ot s_i\in(\fg\ot_k\sqrt{IS}/IS)^\Gamma.$  The ideal $\cI\subseteq (\fg\ot_k S/IS)^\Gamma$ generated by $\alpha$ is contained in $(\fg\ot_k\sqrt{IS}/IS)^\Gamma$, and is thus nilpotent.  The $\cM_i$ are distinct, so by the Chinese remainder theorem, the Lie algebra
$$\cL/(\cM_1\cap\cdots\cap\cM_s)\cong\cL/\cM_1\times\cdots\times\cL/\cM_s$$
is semisimple.  Therefore, $\cI\subseteq\ker\psi$, and the map $\psi$ descends to a surjection 
$$(\fg\ot_k S/IS)^\Gamma/(\fg\ot_k\sqrt{IS}/IS)^\Gamma\rightarrow \cL/\cM_1\times\cdots\times\cL/\cM_s.$$  
Since taking $\Gamma$-invariants is an exact functor, we have
$$(\fg\ot_k S/\sqrt{IS})^\Gamma=(\fg\ot_k S/IS)^\Gamma/(\fg\ot_k\sqrt{IS}/IS)^\Gamma=(\fg\ot_k S/IS)^\Gamma/(\fg\ot_k\sqrt{IS}/IS)^\Gamma,$$
and the map $\psi$ thus gives a surjective Lie algebra homomorphism 
\begin{equation}\label{epimorphism}\fg^M\cong(\fg\ot_k S/M)^{\Gamma^M}\cong(\fg\ot_k S/\sqrt{IS})^\Gamma\stackrel{\psi}{\rightarrow}\cL/\cM_1\times\cdots\times\cL/\cM_s.
\end{equation}\qed

\section{Classification of finite-dimensional simple modules}

We maintain the notation of the previous section: $k$ will be an algebraically closed field of characteristic zero, $S$ a unital, commutative, associative, and reduced $k$-algebra of finite type, $\fg$ a finite-dimensional simple $k$-Lie algebra, and $\Gamma$ a finite group acting by $k$-Lie algebra automorphisms on $\fg\ot_k S$.  By Proposition \ref{1-cocycle}, there then exists an action of $\Gamma$ on $S$ and a $1$-cocycle $u:\ \Gamma\rightarrow \Aut_{S-Lie}(\fg\ot_k S)$ defining the action of $\Gamma$: ${}^\gamma(x\ot s)=u_\gamma(x\ot{}^\gamma s)$, for all $\gamma\in\Gamma$ and $x\ot s\in\fg\ot_k S$.  The action of $\Gamma$ is $R=S^\Gamma$-linear, and the associated twisted current algebra is the $k$-Lie algebra $\cL=(\fg\ot_k S)^\Gamma$. All modules (representations) will be finite dimensional over $k$.  Unless explicitly indicated otherwise, $\ot$ will denote the tensor product $\ot_k$ taken over $k$. 

\begin{theorem}\label{complete-list}
Let $\phi:\ \cL\rightarrow\End_k(V)$ be an irreducible finite-dimensional representation of $\cL$.  Then there exist maximal ideals $M_1,\ldots,M_r\in\Max\,S$ in distinct $\Gamma$-orbits, finite-dimensional simple $\fg^{M_i}$-modules $(V_i,\phi_i)$ (for $i=1,\ldots,r$), and a (possibly trivial) $1$-dimensional $\cL$-module $W$, such that $V\cong W\ot V_1\ot\cdots\ot V_r$, with $\cL$-action defined by $\phi_i\circ \ev_{M_i}$ on each $V_i$.    
\end{theorem}
\proof Since $\cL/\ker\phi$ has a faithful finite-dimensional irreducible representation, it is reductive with centre $Z$ of dimension at most $1$.  See \cite[\S 6 Exercice 20]{bourbaki-Lie1} or \cite[Exercise 12.4]{erdmann-wildon}, for instance.  Hence,
$$\cL/\ker\phi\cong Z\times\cL/\cM_1\times\cdots\times\cL/\cM_\ell,$$
for distinct maximal ideals $\cM_1,\ldots,\cM_\ell\subseteq\cL$ where $\cL/\cM_i$ is a finite-dimensional simple Lie algebra for each $i=1,\ldots,\ell$.  Therefore, $V$ is isomorphic to $W\ot U_1\ot\cdots\ot U_\ell$, where $W$ is a $1$-dimensional $\cL$-module (trivial if $Z$ is zero), and $U_i$ is a finite-dimensional simple $\cL/\cM_i$-module, on which $\cL$ acts by pulling back the action of $\cL/\cM_i$.

Suppose $I(\cM_{i_1})=\cdots=I(\cM_{i_s})$ for some $i_1,\ldots,i_s\in\{1,\ldots,r\}$.  If $M$ is a maximal ideal of $S$ lying over $I(\cM_{i_1})\subset R$, then Proposition \ref{fgM-epi} gives an epimorphism $\psi:\ \fg^M\rightarrow\cL/\cM_{i_1}\times\cdots\times\cL/\cM_{i_s}$.  The kernel of the composition $\psi\circ\ev_{M_i}:\ \cL\rightarrow\cL/\cM_{i_1}\times\cdots\times\cL/\cM_{i_s}$ is contained in $\cM_{i_1}\cap\cdots\cap\cM_{i_s}$, so the irreducible representation
$$\cL\rightarrow\cL/\cM_{i_1}\times\cdots\times\cL/\cM_{i_s}\rightarrow \hbox{End}\,(U_{i_1}\ot\cdots\ot U_{i_s})$$
factors through an irreducible representation of $\fg^M$:
$$\cL\stackrel{\ev_M}{\rightarrow}\fg^M\stackrel{\psi}{\rightarrow}\cL/\cM_{i_1}\times\cdots\times\cL/\cM_{i_s}\rightarrow \hbox{End}\,(U_{i_1}\ot\cdots\ot U_{i_s}).$$

Let $I_1,\ldots,I_r\in\Max\,R$ be the distinct ideals in the set $\{I(\cM_j)\ :\ j=1,\ldots,\ell\}$.  For each $i=1,\ldots,r$, let $M_i$ be any maximal ideal of $S$ lying over $I_i$.  Since the group $\Gamma$ acts on the fibres in $\Max\,S$ over each ideal in $\Max\,R$, the $M_i$ belong to distinct $\Gamma$-orbits.  By the discussion above, the $\cL$-module $V$ is thus isomorphic to $W\ot V_1\ot \cdots\ot V_r$ for finite-dimensional simple $\fg^{M_i}$-modules $V_i=U_{i_1}\ot\cdots\ot U_{i_{n_i}}$, where $M_i$ lies over $I(\cM_{i_j})$ and $U_{i_j}$ is the $\cL/\cM_{i_j}$-module in the factorization $V\cong W\ot U_1\ot\cdots U_\ell$.



\qed

\bigskip

We now prove the converse of Theorem \ref{complete-list}.
\begin{theorem} \label{tensors-are-simple}Let $M_1,\ldots,M_r\in\Max\,S$ be representatives of distinct $\Gamma$-orbits.  Let $(V_i,\phi_i)$ be a finite-dimensional irreducible representation of $\fg^{M_i}$ for $i=1,\ldots,r$, and suppose that $(W,\rho)$ is a 
(possibly trivial) 1-dimensional module for the Lie algebra $\cL$.  Then $V=W\ot V_1\ot\cdots\ot V_r$ is a finite-dimensional simple module for $\cL$ with action $\phi_i\circ\ev_{M_i}$ on each $V_i$.
\end{theorem}
\proof By construction, $V$ is a simple module for $Z\times\fg^{M_1}\times\cdots\times\fg^{M_r}$, where $Z=\ker\rho$.  It thus suffices to show that the map
$$(p,\ev_{M}):\ \cL\rightarrow Z\times\fg^{M_1}\times\cdots\times\fg^{M_r}$$
is surjective, where $p$ is the canonical projection $\cL\rightarrow\cL/\cM_0$ and $\ev_M=(\ev_{M_1},\ldots,\ev_{M_r}).$  Let $I_i=M_i\cap R$ for $i=1,\ldots,r$.  Since $S/R$ is integral and the $M_i$ belong to distinct $\Gamma$-orbits, the $I_i$ are distinct maximal ideals of $R$.  Thus $I_i+I_j=R$ for all $i\neq j$, from which it follows that $I_iS+I_jS=S$ whenever $i\neq j$.  Hence
$$S\left/\bigcap_{i=1}^rI_iS\ \cong\  (S/I_1S)\oplus\cdots\oplus (S/I_rS),\right.$$
and the canonical map 
\begin{equation}
\pi:\ \fg\ot S\ \rightarrow\ \fg\ot S\left/\left(\fg\ot\bigcap_{i=1}^r I_iS\right)\ \right.\cong\ \bigoplus_{i=1}^r\left(\fg\ot (S/I_iS)\right)
\end{equation}
is a surjection.  Taking $\Gamma$-invariants is an exact functor, so the restriction of $\pi$ to $\cL$ is also surjective:
$$\pi:\ \cL=(\fg\ot S)^\Gamma\rightarrow\left(\bigoplus_{i=1}^r(\fg\ot S/I_iS)\right)^\Gamma.$$
Since $I_iS$ is $\Gamma$-stable, $\left(\bigoplus_{i=1}^r(\fg\ot S/I_iS)\right)^\Gamma=\bigoplus_{i=1}^r\left(\fg\ot S/I_iS)\right)^\Gamma$.  By (\ref{reduction-iso}), $(\fg\ot S/I_iS)^\Gamma$ is isomorphic to $\cL/I_i\cL$, so the induced projection 
$$\cL\rightarrow(\cL/I_1\cL)\times\cdots\times(\cL/I_r\cL)$$
is surjective.  But
$$I_i\cL=(M_i\cap R)\cL\subseteq \fg\ot(M_i\cap R)S\subseteq\fg\ot M_i=\ker(\ev_{M_i}),$$
so 
$$\ev_M:\ \cL\rightarrow\ev_{M_1}(\cL)\oplus\cdots\oplus\ev_{M_r}(\cL)=\fg^{M_1}\oplus\cdots\fg^{M_r}$$
is surjective, and the module $V_1\ot\cdots\ot V_r$ is irreducible.  Tensoring with a $1$-dimensional module $W$ does not change irreducibility, so $V=W\ot V_1\ot\cdots\ot V_r$ is also irreducible.

\qed

Theorems \ref{complete-list} and \ref{tensors-are-simple} give a complete list of the finite-dimensional simple modules for any twisted current algebra $\cL$.  However, many of these representations are isomorphic, so we now work at classifying such modules up to isomorphism.  

\bigskip

In Section 2, we introduced an evaluation map $\ev_M$ associated with any maximal ideal $M$ of $S$:
\begin{align*}
\ev_M:\ \Aut_{k-Lie}(\fg\ot S)&\rightarrow\Aut_{k-Lie}(\fg)\\
\psi&\mapsto\psi(M),
\end{align*}
where $\psi(M)x:=(\psi(x\ot 1))(M)$.  Through a small abuse of notation, we associated an automorphism $\gamma(M)$ of $\fg$ with each $\gamma\in\Gamma$ and $M\in\Max\,S$.  Since $\Gamma$ (and thus $\Gamma^M$) is finite, there are finitely many $\gamma(M)$ for any fixed $M$, so the $\Gamma^M$-fixed subalgebra $\fg^M$ is reductive by \cite[\S 1.5 Proposition 14]{bourbaki-Lie7}.  The map $\Gamma\rightarrow\Aut_{k-Lie}(\fg)$ induced by $\ev_M$ is {\em not} a group homomorphism.  However, the group $\Gamma$ acts on $\Aut_{k-Lie}(\fg\ot S)$ as in (\ref{action-on-autos}):
$${}^\gamma\psi=(1\ot\gamma)\circ\psi\circ(1\ot\gamma^{-1}),$$
for $\gamma\in\Gamma$ and $\psi\in\Aut_{k-Lie}(\fg\ot S)$, and with respect to this action, it is easy to verify that the evaluation maps $\ev_M$ give {\em local twisted actions} of $\Gamma$ on $\fg$.  That is, as automorphisms of $\fg$,
\begin{equation}\label{twisted-actions}
\ev_M(\gamma\eta)=\ev_M(\gamma)\circ\ev_M({}^\gamma\eta),
\end{equation}
for all $\gamma,\eta\in\Gamma$ and $M\in\Max\,S$.  It is then straightforward to verify that
\begin{align}
(\gamma(M))^{-1}&=\gamma^{-1}\left({}^{\gamma^{-1}}M\right),\\
\gamma^{-1}(M)&=\left(\gamma({}^\gamma M)\right)^{-1}, \label{3.6}\\
({}^\gamma \eta)({}^\gamma M)&=\eta(M),\label{3.7}
\end{align}
for all $\gamma,\eta$, and $M$.
\begin{lemma}\label{lemma-3.8}
Let $\gamma\in\Gamma$ and $M\in\Max\,S$.  Then $\fg^M=\gamma^{-1}(M)\fg^{{}^\gamma M}$.
\end{lemma}
\proof Let $\mu\in\Gamma^{{}^\gamma M}$ and $x\in\fg^M$.  Then $\gamma^{-1}\mu\gamma\in\Gamma^M$, so $x=(\gamma^{-1}\mu\gamma)(M)x$.  By (\ref{twisted-actions}), $(\gamma^{-1}\mu\gamma)(M)=\gamma^{-1}(M){}^{\gamma^{-1}}\mu(M){}^{\gamma^{-1}\mu}\gamma(M)$, so by (\ref{3.7}) we have
\begin{align*}
(\gamma^{-1}(M))^{-1}x&={}^{\gamma^{-1}}\mu(M){}^{\gamma^{-1}\mu}\gamma(M)x\\
&=\mu({}^\gamma M)\gamma({}^{\mu^{-1}\gamma}M)x.
\end{align*}  
But $\mu\in\Gamma^{{}^\gamma M}$, so $\mu^{-1}\in\Gamma^{{}^\gamma M}$, and thus 
\begin{align*}
(\gamma^{-1}(M))^{-1}x&=\mu({}^\gamma M)\gamma({}^\gamma M)x\\
&=\mu({}^\gamma M)[\gamma^{-1}(M)]^{-1}x,
\end{align*}
by (\ref{3.6}).  Therefore, $(\gamma^{-1}(M))^{-1}x\in\fg^{{}^\gamma M}$ and 
\begin{equation}\label{3.9}
\fg^M\subseteq \gamma^{-1}(M)\fg^{{}^\gamma M}.
\end{equation}
Replacing $\gamma$ by $\gamma^{-1}$, and then $M$ by ${}^\gamma M$ in (\ref{3.9}), we obtain: $\fg^{{}^\gamma M}\subseteq\gamma({}^\gamma M)\fg^M.$  The reverse inclusion to (\ref{3.9}) then follows from (\ref{3.6}).
\qed

\bigskip

The finite-dimensional simple $\cL$-modules are precisely the tensor products of $1$-dimensional $\cL$-modules with irreducible representations of finite-dimensional semisimple quotients of $\cL$.  In Theorem \ref{complete-list}, we saw that the latter factor through evaluation maps $\cL\stackrel{\ev_{\underline{M}}}{\longrightarrow}\fg^{M_1}\oplus\cdots\oplus\fg^{M_r}$, for some maximal ideals $M_1,\ldots, M_r$ of $S$ from distinct $\Gamma$-orbits.  Pullbacks of $\fg^{M_1}\oplus\cdots\oplus\fg^{M_r}$-modules are called {\em evaluation representations} of $\cL$.  To classify the finite-dimensional simple $\cL$-modules up to isomorphism, it thus suffices to classify irreducible evaluation representations up to isomorphism and determine which ones appear as pullbacks of representations of semisimple quotients of $\cL$.

Let $\Rep(\fg^M)$ be the set of isomorphism classes of finite-dimensional irreducible representations of the (reductive) Lie algebra $\fg^M$.  Consider the fibre bundle 
$$\mathcal{B}=\coprod_{M\in\Max\,S}\Rep(\fg^M)\rightarrow\Max\,S,$$
where the fibre over $M\in\Max\,S$ is $\Rep(\fg^M)$.  The class $[\mathbbm{1}]$ of the trivial $1$-dimensional module $\mathbbm{1}$ is a member of each fibre.  For any section $\psi:\ \Max\,S\rightarrow\mathcal{B}$ of this bundle, we define the {\em support} of $\psi$ to be the set
$$\supp\psi=\{M\in\Max\,S\ :\ \psi(M)\neq [\mathbbm{1}]\}.$$
A section $\psi$ is {\em finitely supported} if $\supp\psi$ is of finite cardinality.

By Lemma \ref{lemma-3.8}, $(\gamma(M))^{-1}\fg^M=\fg^{{}^{\gamma^{-1}}M}$ for all $\gamma\in\Gamma$ and $M\in\Max\,S$.  We write $\psi({}^{\gamma^{-1}}M)\circ(\gamma(M))^{-1}$ for the isomorphism class $[\phi\circ(\gamma(M))^{-1}]\in\Rep(\fg^M)$, where $\phi$ is any representative of the isomorphism class $\psi({}^{\gamma^{-1}}M)\in\Rep(\fg^{{}^{\gamma^{-1}}M})$.

\begin{lemma}
The group $\Gamma$ acts on the set $\mathcal{F}$ of finitely supported sections of $\mathcal{B}$ by
\begin{equation}\label{action}
(\gamma\cdot\psi)(M)=\psi({}^{\gamma^{-1}}M)\circ(\gamma(M))^{-1},
\end{equation}
for all $\gamma\in\Gamma$, $M\in\Max\,S$, and $\psi\in\mathcal{F}$.
\end{lemma}
\proof This is a straightforward verification.\qed

\bigskip

Let $\psi$ be an element of $\mathcal{F}^\Gamma$, the set of finitely supported sections fixed under the action (\ref{action}):
$$\psi\in\mathcal{F}^\Gamma=\{\psi\in\mathcal{F}\ :\ \gamma\cdot\psi=\psi\hbox{\ for all\ }\gamma\in\Gamma\}.$$
We now establish a bijection $T$ between $\mathcal{F}^\Gamma$ and the set $\mathcal{S}$ of isomorphism classes of irreducible evaluation representations of $\cL$.  For each $\psi\in\mathcal{F}^\Gamma$, we see that $\supp\psi$ is a union of (finitely many) distinct $\Gamma$-orbits $\Omega_1,\ldots ,\Omega_r$ in $\Max\,S$.  Choose $\{M_1,\ldots,M_r\}$ a set of representatives of these orbits.  Consider the isomorphism class $[\rho_\psi]$ of the representation
$$\rho_\psi:\ \cL\stackrel{\ev_{\underline{M}}}{\longrightarrow}\fg^{M_1}\oplus\cdots\oplus\fg^{M_r}\stackrel{(\rho_1,\ldots,\rho_r)}{\longrightarrow}\End\,(V_1\ot\cdots\ot V_r),$$
where $(V_i,\rho_i)$ is a representative of the isomorphism class $\psi(M_i)\in\Rep(\fg^{M_i})$.
\begin{lemma} For every $\psi\in\mathcal{F}^\Gamma$, the isomorphism class $[\rho_\psi]$ depends only on $\psi$ and is independent of the choice of representatives $M_i$ and $(V_i,\rho_i)$ for all $i$.  That is, there is a well-defined map $T:\ \mathcal{F}^\Gamma\rightarrow\mathcal{S}$ given by $T(\psi)=[\rho_\psi].$
\end{lemma}
\proof Suppose $M_i'$ is a representative of the same $\Gamma$-orbit as $M_i$ for some $i$.  Then $M_i'={}^\gamma M_i$ for some $\gamma\in\Gamma$, and $\rho_i\circ\ev_{M_i}:\ \cL\rightarrow\End\,V_i$ and $\rho_i'\circ\ev_{M_i'}:\ \cL\rightarrow\End\,V_i$ determine exactly the same representation, where $\rho_i'=\rho_i\circ\gamma^{-1}(M_i)$.  Indeed, for all $z=\sum_j x_j\ot s_j\in\cL$, we have
\begin{align*}
(\gamma^{-1}(M_i)\circ\ev_{M_i'})(z)&=\gamma^{-1}(M_i)\sum_j s_j(M_i')x_j\nonumber\\
&=\gamma^{-1}(M_i)\sum_j{}^{\gamma^{-1}}s_j(M_i)x_j\quad\quad\ \, \quad\quad\quad\quad\quad\quad\quad\quad\quad\quad\hbox{ (by (\ref{1star}))}\nonumber\\
&=\sum_j{}^{\gamma^{-1}}s_j(M_i){}^{\gamma^{-1}}(x_j\ot 1)(M_i)\nonumber\\
&=\sum_j\left({}^{\gamma^{-1}}s_ju_{\gamma^{-1}}(x_j\ot 1)\right)(M_i)\nonumber\\
&=\sum_ju_{\gamma^{-1}}(x_j\ot{}^{\gamma^{-1}}s_j)(M_i)\quad\quad\ \, \quad\quad\quad\quad\quad\quad\quad\quad\quad\quad\hbox{ (by (\ref{2stars}))}\nonumber\\
&={}^{\gamma^{-1}}\left(\sum_j x_j\ot s_j\right)(M_i)\\
&={}^{\gamma^{-1}}z(M_i)\\
&=z(M_i)\quad\quad\ \, \quad\quad\quad\quad\quad\quad\quad\ \quad\quad\quad\quad\quad(\hbox{since\ }z\in\cL=(\fg\ot S)^\Gamma)\nonumber\\
&=\ev_{M_i}(z).
\end{align*}
Hence, as maps on $\cL$,
\begin{align*}
\rho_i\circ\ev_{M_i}&=\rho_i\circ\gamma^{-1}(M_i)\circ\ev_{M_i'}\\
&=\rho_i'\circ\ev_{M_i'},
\end{align*}
and $T$ is well-defined.\qed

\bigskip

\begin{proposition}\label{equimaps}
The map $T:\ \mathcal{F}^\Gamma\rightarrow\mathcal{S}$ is a bijection.
\end{proposition}
\proof By Theorem \ref{complete-list}, $T$ is surjective.  It is straightforward to verify that $T$ is also injective, using the argument in \cite[Proposition 4.15]{NSS}.\qed

\begin{theorem}
The isomorphism classes of finite-dimensional simple $\cL$-modules are in natural bijection with the pairs $(\lambda,\psi)\in\cL^*\times\mathcal{F}^\Gamma$, such that $\lambda$ vanishes on $[\cL,\cL]$, $\cL/\ker(\rho_\psi)$ is semisimple, and $\ker(\lambda+\rho_\psi)=\ker(\lambda)\cap\ker(\rho_\psi)$ for any $\rho_\psi\in T(\psi)$.
\end{theorem}
\proof By the argument in Theorem \ref{complete-list}, every finite-dimensional simple $\cL$-module $(U,\phi)$ is the pullback of a representation $W\ot V$ of a (reductive) Lie algebra $\ol{\cL}=\cL/\ker\phi=Z\oplus L$, where $Z$ is the centre and $L$ is the derived subalgebra of $\ol{\cL}$.  Here $W$ is a $1$-dimensional $Z$-module, on which $Z$ acts by scalars determined by a linear functional $\lambda\in Z^*$.  Since $Z=(Z\oplus L)/ L=\ol{\cL}/[\ol{\cL},\ol{\cL}]\cong\cL/([\cL,\cL]+\ker\phi)$, such $\lambda$ may be identified with linear functionals on $\cL$, which vanish on $[\cL,\cL]+\ker\phi$.  

By Proposition \ref{equimaps}, the $\cL$-module $(V,\rho)$ obtained by pulling back the $L$-action on $V$ to $\cL$ is then described, up to isomorphism, by a $\Gamma$-invariant map $\psi\in\mathcal{F}^\Gamma$.  It is clear that $\ker\rho$ is an ideal of $\cL$ which contains both $\ker\phi$ and the preimage $\mathcal{Z}\subseteq\cL$ of $Z$ under the projection $\cL\rightarrow\ol{\cL}$.  That is, $\cL/\ker\rho$ is a quotient of $L$, a semisimple Lie algebra, so $\cL/\ker\rho$ is semisimple.  Moreover, $W$ is $1$-dimensional, so $\ker\phi=\ker(\lambda+\rho)$.  Thus $\ker(\lambda)\supseteq\ker(\lambda+\rho)$, from which it follows that $\ker(\lambda+\rho)=\ker(\lambda)\cap\ker(\rho)$.  Conversely, any pair $(\lambda,\psi)$ satisfying the conditions of the theorem will determine a finite-dimensional irreducible representation $(W\ot V,\phi)$, where $\lambda$ determines a representation $W$ of the centre of $\overline{\cL}=\cL/\ker\phi$ and $\psi$ gives a representation $V$ of its derived subalgebra $[\overline{\cL},\overline{\cL}]$.

Suppose the representations $(W\ot V,\phi)$ and $(W'\ot V',\phi')$ determined by two such pairs $(\lambda,\psi)$ and $(\lambda',\psi')$ are isomorphic $\cL$-modules.  Then $\ker\phi=\ker\phi'$ and the factorizations $W\ot V$ and $W'\ot V'$ correspond to the same (unique) decomposition $\ol{\cL}=\cL/\ker\phi=Z\oplus L$ described above.  The representations $W\ot V$ and $W'\ot V'$ are isomorphic as $\ol{\cL}$-modules, and by \cite[\S 7.7 Proposition 8]{bourbaki-algebre8}, $W\cong W'$ as $Z$-modules and $V\cong V'$ as $L$-modules.  By the discussion above, we then see that $\lambda=\lambda'$.  Similarly, the pullbacks of $V$ and $V'$ to $\cL$ are isomorphic $\cL$-modules, so $\psi=\psi'$ by Proposition \ref{equimaps}.\qed

\end{document}